\documentclass{amsart}

\usepackage{amsmath,amssymb}

\newtheorem{theorem}{Theorem}
\newcommand{\bt}{\begin{theorem}}
\newcommand{\et}{\end{theorem}}
\newtheorem{lemma}{Lemma}
\newcommand{\bl}{\begin{lemma}}
\newcommand{\el}{\end{lemma}}
\newtheorem{corollary}{Corollary}
\newcommand{\bc}{\begin{corollary}}
\newcommand{\ec}{\end{corollary}}
\newcommand{\beq}{\begin{equation}}
\newcommand{\eeq}{\end{equation}}
\newcommand{\benum}{\begin{enumerate}}
\newcommand{\eenum}{\end{enumerate}}

\newcommand{\Z}{\ensuremath{\mathbf Z}}
\newcommand{\Q}{\ensuremath{\mathbf Q}}
\newcommand{\R}{\ensuremath{\mathbf R}}

\newcommand{\mbe}{\ensuremath{ \mathbf e}}
\newcommand{\mbv}{\ensuremath{ \mathbf v}}

\DeclareMathOperator{\length}{length}

\newcommand{\bmat}{\left(\begin{matrix}}
\newcommand{\emat}{\end{matrix}\right)}

\DeclareMathOperator{\qqand}{\qquad\text{and}\qquad}

\title[Decomposition of parallelepipeds]{Von Neumann-Steinhaus decomposition of parallelepipeds}

\author{Melvyn B. Nathanson}

\address{Department of Mathematics\\Lehman College (CUNY)\\Bronx, NY 10468} 

\email{melvyn.nathanson@lehman.cuny.edu}



\thanks{Supported in part by a grant from the PSC-CUNY Research Award Program.}

\date{\today}
\begin{document}

\begin{abstract}
A problem of Steinhaus was to partition a finite interval $I$ of the real line 
into countably infinitely many pairwise disjoint 
sets that are congruent in the sense that each set is a translate 
of a fixed set $A$.  This paper describes von Neumann's solution of the Steinhaus problem, 
and an extension to parallelepipeds in $\R^n$.  
\end{abstract}

\maketitle

\section{ A problem of Hugo Steinhaus}

In a beautiful essay, ``A walk through Johnny von Neumann's garden,'' Freeman Dyson~\cite{dyso13} 
describes some of von Neumann's most important contributions to mathematics, physics, and computer science.  He writes:    
\begin{quotation}
In another corner of the garden, there is a little flower all by itself, a short paper 
\ldots [that] solves a problem raised by the Polish mathematician Hugo Steinhaus\ldots.  
Johnny solved the Steinhaus problem quickly and never returned to it.  
The theorem that he proved is counterintuitive, and the proof is astonishing. 
\end{quotation}
Steinhaus's problem (first stated in 1921 in a paper of Stefan Mazurkiewicz~\cite{mazu21}) 
was to partition a finite interval of the real line 
into countably infinitely many pairwise disjoint 
sets that are congruent in the sense that each set is a translate 
of a fixed set $A$.  

In 1928, von Neumann~\cite{vonn28} solved Steinhaus's problem 
for open, for closed, and for half-open, half closed intervals. 
For every finite interval $I$, von Neumann proved the existence 
of a set $A$ and a sequence of real numbers $(b_n)_{n=1}^{\infty}$ 
such that  the congruent sets 
\[
\{A+b_n: n = 1,2,3,\ldots\} 
\]
are pairwise disjoint and 
\[
I = \bigcup_{n=1}^{\infty} (A + b_n).
\]   
The set $A$  is necessarily non-measurable.   

This note describes von Neumann's almost forgotten ``little flower.''
It gives von Neumann's proof for intervals in \R\ 
(perhaps, the first exposition in English) 
and extends the result to parallelepipeds in $\R^n$.

\section{Additive number theory}

Steinhaus's problem is easily and naturally restated in the language of additive number theory. 
Let $\Omega$ be a semigroup, written additively, 
and let $A$, $B$, and $C$ be subsets of  $\Omega$.  
We define the \emph{sumset} 
\[
A+B = \{a+b: a \in A \text{ and } b \in B\}.
\] 
The set $C$ is the \emph{sum} of the sets $A$ and $B$ if 
\[
A+B=C.
\]
We have the additive identity:
\[
C = A+B= \bigcup_{a\in A} (a+B)  = \bigcup_{b\in B} (A+b) .
\]
The translates $a+B$ are congruent for all $a \in A$ 
and the  translates  $A+b$ are congruent for all $b \in B$.   
It follows that the set $C$ is the union of countably infinitely 
many congruent sets of the form  $A+b$
if and only if $C$ is also the union of any number of  congruent countably infinite sets
of the form $a+B$.

The set $C$ is the \emph{direct sum} of the sets $A$ and $B$, denoted 
\[
A \oplus B = C
\]
if $A+B = C$ and if every element of $C$ has a \emph{unique} 
representation in the form $a+b$ 
with $a \in A$ and $b \in B$.
Thus, $A \oplus B = C$ means that if $(a, b) \in A \times B$ and $(a', b') \in A \times B$  satisfy 
\[
a +b = a'+ b' 
\]
then $(a, b) = (a', b')$. 
We have $A \oplus B = C$ if and only if the sets $\{A+b:b \in B\}$ are pairwise disjoint 
 if and only if the sets $\{a + B: a \in A\}$ are pairwise disjoint.  

This remark gives an equivalent form of Steinhaus's problem.
The original problem: 
Prove that every interval is the union of a countably infinite number 
of pairwise disjoint  uncountable congruent sets.  
The equivalent problem: Prove that every interval is the union of uncountably many pairwise disjoint sets 
that are congruent to a countably infinite set.  

The following elementary results are used in Section~\ref{vonN:section:parallelepiped}. 

\bl                                       \label{vonN:lemma:k-sum} 
Let $\Omega_1, \ldots, \Omega_n$ be semigroups, written additively, 
and let $\Omega_1 \times \cdots \times \Omega_n$ be the semigroup of $n$-tuples   
 with addition defined coordinate-wise.   
Let $A_i$, $B_i$, and $C_i$ be nonempty subsets of $\Omega_i$ for $i \in \{1,2,\ldots, n\}$ 
and let $A_1\times \cdots \times A_n$, $B_1\times \cdots \times B_n$, 
and $C_1\times \cdots \times C_n$ be the corresponding subsets 
of $\Omega_1 \times \cdots \times \Omega_n$.   
Then  
\[
A_i \oplus B_i = C_i
\]
for all $i \in \{1,\ldots, n\}$ if and only if 
\[              
(A_1\times \cdots \times A_n) \oplus (B_1\times \cdots \times B_n) = C_1\times \cdots \times C_n. 
\] 
\el

\bl                                             \label{vonN:lemma:Tisom} 
Let $\Omega$ and $\Omega'$ be semigroups, written additively, 
and let $T:\Omega \rightarrow \Omega'$ be a one-to-one semigroup homomorphism.  
Let $A$, $B$, and $C$ be  nonempty subsets of $\Omega$.  
Define the subsets $A' = T(A)$, $B'=T(B)$, and $C'=T(C)$ of $\Omega'$.  
Then  $A\oplus B = C$ if and only if $A' \oplus B' = C'$. 
\el

\section{Construction of the set $B$ and the group $G$}
Von Neumann proved the following theorem. 

\bt                    \label{vonN:theorem:Steinhaus}  
Every finite interval $I$ is the union of countably infinitely many 
pairwise disjoint congruent sets. 
Equivalently, for every finite interval $I$ there is an uncountable set $A$ 
and a countably infinite set $B$ 
such that $A \oplus B = I$. 
\et

Let $\Omega$ be  a $\sigma$-algebra of subsets of \R\ that contains the intervals of \R.  
A \emph{normalized measure}\index{normalized measure} on $\Omega$ 
is an extended real-valued countably additive nonnegative function $\mu$ on $\Omega$ 
such that the measure of an interval is the length of the interval.  
The measure $\mu$ is translation-invariant if $\mu(A+b) = \mu(A)$ 
for all $A \in \Omega$ and $b \in \R$. 
In 1905 Giuseppe Vitali~\cite{vita05} (cf. von Neumann~\cite[Chapter VII]{vonn50}) 
proved the existence of subsets of 
 the additive abelian group \R\ 
that are not measurable with respect to any normalized translation-invariant 
measure on $\Omega$.  
Vitali's proof uses the axiom of choice to obtain  
a set of coset representatives of the quotient group $\R/\Q$ in the unit interval $[0,1]$.  
Note that \Q\ is a countably infinite and dense subgroup of \R. 

Adapting Vitali's idea, von Neumann applied the axiom of choice 
to the intersections of an interval $I$ 
with the cosets of a different countably infinite and dense subgroup $G$ of  \R.  
The key to the proof is the  construction, for every $\varepsilon > 0$, 
of a countably infinite set $B$ 
of \Q-independent real numbers with exactly two accumulation points.
The subgroup $G$ is the group generated by $B$, that is, 
the set of all integral linear combinations of elements of $B$.

\bl                                   \label{vonN:lemma:Bset} 
Let  $\beta_0$ and $\beta_1$ be real numbers with $\beta_0 < \beta_1$ 
and let 
\[
0 < \varepsilon < \frac{\beta_1 - \beta_0}{2}. 
\] 
There exists a countably infinite set 
of real numbers 
$B = \{b_0,b_1,b_2, b_3, \ldots\}$ such that 
\benum 
\item
the set $B$ is linearly independent over \Q, 
\item 
for all $n =  1, 2, 3, \ldots$, 
\[
\beta_0 < b_0 < b_{2n} < \beta_0 + \varepsilon 
\]
and 
\[
 \beta_1 - \varepsilon <  b_{2n+1}   < b_1 < \beta_1
\]
\item
 \[
 \lim_{n\rightarrow \infty} b_{2n} = b_0 
\qqand 
 \lim_{n\rightarrow \infty} b_{2n+1} = b_1 
 \]
 \item
for all $\delta > 0$, the set 
\[
\{b_n \in B: b_0+\delta < b_n < b_1 - \delta \}
\]
is finite,
  \item
the additive abelian subgroup $G$ of $\R$ that is generated by $B$ 
is countably infinite and dense in \R.  
\eenum
\el

\begin{proof}
Let $ \{c_0, c_1, c_2, \ldots\}$ be any countably infinite set of real numbers that are 
linearly independent over \Q.  
For example, if $\theta = \pi$ or $e$ or any transcendental real number, then the set 
$\{1, \theta, \theta^2, \theta^3, \ldots\}$ is linearly independent over \Q.  
If $(q_n)_{n=0}^{\infty}$ is any sequence of nonzero rational numbers, 
then the set $ \{q_0c_0, q_1c_1, q_2c_2, \ldots\}$ is also linearly independent over \Q.

Choose nonzero rational numbers $q_0$ and $q_1$ so that the numbers 
\[
b_0 = q_0 c_0 \qqand b_1 = q_1 c_1
\]
 satisfy the inequalities 
\[
\beta_0 <  b_0    < \beta_0 + \varepsilon \qqand  \beta_1 - \varepsilon <  b_1  < \beta_1. 
\]
Choose $\varepsilon' > 0$ such that 
\[
  b_0 < b_0 + \varepsilon'    < \beta_0 + \varepsilon \qqand  \beta_1 - \varepsilon <  b_1 - \varepsilon' < b_1. 
\]
For every positive integer $n$, choose nonzero rational numbers $q_{2n}$ and $q_{2n+1}$ 
such that the numbers 
\[
b_{2n} = q_{2n}c_{2n}\qqand b_{2n+1} = q_{2n+1} c_{2n+1} 
\]
satisfy  the inequalities 
\[
b_0 < b_{2n} < b_0 + \frac{\varepsilon'}{n} \qqand b_1 - \frac{\varepsilon'}{n} < b_{2n+1} < b_1.  
\]
The sequence $(b_n)_{n=0}^{\infty}$ satisfies properties (1), (2), and (3). 
Properties (2) and (3) imply property (4).

The subgroup $G$ of \R\ generated by $B$ is the set of all integral 
linear combinations of elements of the countably infinite set $B$ 
and so $G$ is countably infinite.
The number $\theta = b_1/b_0$ is irrational because $b_0$ and $b_1$ are \Q-independent.  
By Kronecker's theorem (Hardy and Wright~\cite[Theorem 439]{hard-wrig}), 
the set of numbers $\{n_0+n_1\theta: n_0,n_1\in \Z\}$ is dense in \R\ and so 
the set of numbers 
\[
\{b_0(n_0+n_1\theta) : n_0,n_1 \in \Z\} = \{n_0 b_0+n_1b_1: n_0,n_1 \in \Z\}
\]
 is also dense in \R.  
The group $G$ contains this set and so $G$ is dense in \R. 
This proves~(5).   
\end{proof}

The following result, implicit in~\cite{vonn28}, is used in the proof of Theorem~\ref{vonN:theorem:Steinhaus}. 

\bl                                                                      \label{vonN:lemma:B-cover} 
Let $\varepsilon > 0$ and let $\beta_0 = -2\varepsilon$ and $\beta_1 =2 \varepsilon$.  
Let 
\[
B = B(\varepsilon) = \{b_0,b_1,b_2, b_3, \ldots\} 
\] 
be the set constructed in Lemma~\ref{vonN:lemma:Bset}, and let $G = G(\varepsilon)$ 
be the additive subgroup of \R\ generated by $B$.  
Let $J$ be a finite interval of $\length(J) > 8 \varepsilon$. 
No finite union of pairwise disjoint translates of $B$ contains the set $J \cap G$.  
\el

\begin{proof}
We have $b_0 = \min(B)$ and $b_1 = \max(B)$, and 
\[
2\varepsilon < b_1 - b_0 < 4\varepsilon. 
\]  
By Lemma~\ref{vonN:lemma:Bset}, 
the set $B$ has the property that, for all $\delta > 0$, the set 
\[
B \cap \left( [b_0, b_0 + \delta) \cup (b_1-\delta, b_1] \right) 
\]
is a cofinite subset of $B$, and so, for all $a \in \R$ 
and $\delta > 0$, the set   
\[
[a + b_0, a + b_0 + \delta) \cup (a + b_1-\delta, a + b_1]
\]
is  a cofinite subset of the translate $a + B$.  Moreover,
\[
\{ a+b_0, a+b_1 \} \subseteq a+B \subseteq [a+b_0, a+b_1]. 
\]

Let $a_1, \ldots, a_k$ be real numbers such that the $k$ translates   
$a_1+B,\ldots, a_k+B$ are pairwise disjoint.    
We shall prove that 
\[
X = \bigcup_{i=1}^k (a_i+B)
\]
does not contain $J \cap G$. 
The endpoints of these translates, that is, the 
$2k$ numbers $a_i + b_j$ for $i \in \{1,\ldots, k\}$ and $j \in \{0,1\}$, are distinct.  
Choose $\delta$ such that 
\[
0 < 3 \delta <  \min \left\{ \left| \left( a_{i}+b_{j}\right) -  \left(a_{i'}+b_{j'} \right) \right| : 
(i,j) \neq (i', j') \right\}.
\]
This inequality implies that 
the $2k$ intervals $[a_i + b_0, a_i + b_0 + \delta)$ 
and $(a_i + b_1 -  \delta, a_i + b_1]$ for $i  \in \{1,\ldots, k\}$ 
are also  pairwise disjoint.  
Moreover, consecutive intervals are separated by intervals of 
length at least $\delta$.  

We have 
\[
0 < 3\delta < b_1-b_0 < \beta_1-\beta_0 = 4\varepsilon < \length(J)  
\]
and so the interval $J$ contains a subinterval $J^*$ 
of length at least $\delta$ that is disjoint from the set 
\[
X^* =  \bigcup_{i=1}^k [a_i + b_0 , a_i + b_0 + \delta) \cup 
\bigcup_{i=1}^k (a_i + b_1 - \delta , a_i + b_1]. 
\] 
 If $X$ contains $J\cap G$, then $X\setminus X^*$ contains $J^* \cap G$.  
Because $G$ is dense in \R, the interval $J^*$ contains infinitely many points of $G$, 
that is, $J^* \cap G$ is infinite. 
However, the set $X^*$ contains all but finitely many elements of $X$, that is, 
$X\setminus X^*$ is finite.  
The finite set $X\setminus X^*$ does not contain the infinite set $J^* \cap G$, 
and so $X$ does not contain $J\cap G$. This completes the proof. 
\end{proof}

\section{Beginning the proof of von Neumann's theorem}

Let $\varepsilon > 0$ and let 
$\beta_0 = -2\varepsilon$ and $\beta_1 = 2\varepsilon$. 
Apply Lemma~\ref{vonN:lemma:Bset} to obtain the set 
\[
B = B(\varepsilon) = \{b_0, b_1, b_2,\ldots \} \subseteq (-2\varepsilon, 2\varepsilon) 
\]
with accumulation points $b_0 \in (-2\varepsilon, -\varepsilon)$ 
and $b_1 \in (\varepsilon, 2\varepsilon)$. 
Let  $G = G(\varepsilon)$ be the countably infinite and dense subgroup of \R\ generated by $B$. 
Apply the axiom of choice to obtain  
an uncountable set $V$ of coset representatives of the quotient group $\R/G$.  
Then  
\[
\R = \bigcup_{v \in V} (v+G). 
\]
 For all $v \in V$, the coset $v+G$ in the quotient group $\R/G$ 
 is also countably infinite and dense in \R.
 
 Let $I$ be an open, closed, or half-open, half-closed interval  with 
 \[
0 < 8\varepsilon < \length(I). 
 \]  
For all $v \in V$ the translate 
\[
I_v = I - v = \{x-v :x\in I\}
\]
 is an interval in \R\ with $\length(I_v) = \length(I)$.  
We have $I = v + I_v$ and 
\[
 I \cap (v+G) =  (v+I_v) \cap (v+G) = v+ ( I_v \cap G).
\]
Because the cosets $v+G \in \R/G$ partition \R, the sets $\{I\cap (v+G): v \in V\}$ 
partition $I$ and 
\[
I  = \bigcup_{v \in V}   \left(  I \cap (v+G)  \right)
 = \bigcup_{v \in V} \left(  v+ \left( I_v \cap G  \right) \right). 
\]
Suppose that for all $v \in V$ there is a set $A_v$ such that 
\[
I_v \cap G = A_v \oplus B. 
\]
Let 
\[
A = \bigcup_{v\in V} (v+A_v). 
\]
We obtain 
\begin{align*}
I & = \bigcup_{v \in V} \left(  v+ \left( I_v \cap G  \right) \right) \\ 
& = \bigcup_{v \in V} \left(  v+ \left(  A_v \oplus B  \right) \right) \\ 
& = \bigcup_{v \in V} \left(  \left(  v+  A_v   \right) \oplus  B  \right) \\
& =  \left(   \bigcup_{v \in V} \left(  v+  A_v  \right)   \right) \oplus  B \\
& = A \oplus B.
\end{align*}
The sumset $A\oplus B = I$ is a direct sum because 
$\left(  v+  A_v   \right) \oplus  B = I\cap (v+G)$  
and the sets $\{I\cap (v+G): v \in V\}$ 
are pairwise disjoint. 
Thus, to solve Steinhaus's problem, it suffices to prove the following result. 

\bt               \label{vonN:theorem:JG=AB}
Let $J$ be a finite interval.  
Let 
\[
0 < 8 \varepsilon < \length(J)   
\]
and let $B = B(\varepsilon)$ be the set and $G = G(\varepsilon)$ the group constructed in Lemma~\ref{vonN:lemma:Bset} from 
$\beta_0 = -2\varepsilon$ and $\beta_1 = 2\varepsilon$.   
Let $C = J \cap G$.  
There is a set $A$ such that 
\[
 A \oplus B = C. 
\]
\et

\section{A core lemma} 

The proof of von Neumann's theorem is inductive and based on the following lemma. 

\bl                  \label{vonN:lemma:core}
Let $c_0$, $c_1$, and $\varepsilon$ be real numbers with 
\[
0 < \varepsilon <  \frac{c_1-c_0}{8}.
\]
Let $B = B(\varepsilon)$ be the set obtained from 
Lemma~\ref{vonN:lemma:Bset} with $\beta_0 = -2\varepsilon$ 
and $\beta_1 = 2\varepsilon$, and let $G = G(\varepsilon)$ 
be the additive abelian group generated by $B$.  
Let $J'$ be the open interval $(c_0, c_1)$. 
Let
\[
x \in J' \cap G   
\]  
and let $a_1,\ldots, a_k$ be real numbers such that 
\[
x \notin \bigcup_{i=1}^k (a_i + B). 
\]
There exists a  real  number $a$ such that 
\beq            \label{vonN:aB11}
x \in a  + B \subseteq J' \cap G
\eeq
and
\beq            \label{vonN:aB22}
\left( a + B  \right) \cap  \bigcup_{i=1}^k  (a_i + B)  = \emptyset.
\eeq
\el

\begin{proof}
We shall prove that property~\eqref{vonN:aB11} holds for infinitely many 
numbers $a $ and that property~\eqref{vonN:aB22} 
excludes only finitely many numbers $a$ for which property~\eqref{vonN:aB11} holds. 

Let  $B = \{b_n: n= 0,1,2,\ldots\}$, where 
 $b_0 = \min(B)$ and $b_1 = \max(B)$ are the accumulation points of $B$, 
and, for all positive integers $n$,  there are the inequalities 
\[
-2\varepsilon < b_0 < b_{2n} < -\varepsilon 
\]
and 
\[
 \varepsilon < b_{2n+1} < b_1 < 2\varepsilon.  
\]
It follows that  $2\varepsilon < b_1 - b_0 < 4\varepsilon$. 

We have $x \in a+B$ if and only if $x = a+b_j$ for some $b_j \in B$ 
if and only if $a=x-b_j$ for some $b_j \in B$.  
Because $x \in G$ and the set $B$ generates $G$, it follows that, 
for all $b_j \in B$, if $a=x-b_j$, then  
\[
a+B = \{a +b_n: b_n \in B\} = \{x-b_j+b_n: b_n \in B\} \subseteq G.  
\]
We shall prove  that $a+B \subseteq J'$ for infinitely many $a$, or, equivalently, 
that  
\[
\{x-b_j+b_n: b_n \in B\} \subseteq J' 
\]
for infinitely many $b_j \in B$.  
For all $b_n \in B$ we have 
\[
x - b_j + b_0 \leq x-b_j + b_n   \leq x-b_j + b_1.    
\]
Therefore, $ x - b_j + B \subseteq (c_0, c_1)$ 
if 
\begin{align*}
c_0  < x - b_j + b_0 \qqand x-b_j + b_1 < c_1  
\end{align*} 
or, equivalently, if 
\[
x - c_1 + b_1 < b_j < x - c_0 + b_0.    
\]   
Thus, it suffices to prove that the open interval 
\[
J'' = \left( x - c_1 + b_1, \ x - c_0 + b_0  \right) 
\]
 contains infinitely many numbers $b_j \in B$. 
 We shall prove that $J''$ contains either $b_0$ or $b_1$. 
 Because $b_0$ and $b_1$ are accumulation points of $B$, it follows 
 that  the open interval $J''$ contains infinitely many numbers $b_j \in B$.  

The interval $J''$ has length  
\begin{align*}
\length(J'') & = (x - c_0 + b_0) - ( x - c_1 + b_1 ) \\
& = (c_1 - c_0)  - (b_1 - b_0 )  
 > 8\varepsilon - 4 \varepsilon \\ 
& = 4\varepsilon  >b_1 - b_0. 
\end{align*}
Because  $x \in J' =  (c_0, c_1)$, we have 
\[
x - c_1 + b_1 <  b_1 
\qqand 
b_0 < x-c_0+b_0. 
\]
There are two cases.  In the first case,  $  x-c_1+b_1 < b_0$ and 
\[
x-c_1+b_1 < b_0 < x-c_0+b_0 
\]
and so $b_0 \in J''$.

In the second case,  $x-c_1+b_1 \geq b_0$.  If $b_1 \geq x-c_0+b_0$, then $J'' \subseteq (b_0, b_1)$ 
and $b_1 - b_0 < \length(J'') \leq  b_1-b_0$, which is absurd.  Therefore, 
\[
x - c_1 + b_1  < b_1 < x-c_0+b_0 
\]
and $b_1 \in J''$.  
This proves that  $x \in a + B \subseteq J' \cap G$ for infinitely many $a \in \R$.  

Next we prove that there are only finitely many $a$ such that 
\[
x \in a+B \qqand (a+B) \cap \bigcup_{i=1}^k (a_i+B) \neq \emptyset. 
\]
For all $i \in \{1,\ldots, k\}$, the set $a_i+B$ is a translate of $B$ that does not contain $x$.
We shall prove that there are at most two translates  of $B$ 
that both contain $x$ and intersect $a_i+B$. 

Let $i \in \{1,2,\ldots, k\}$.  For $h \in \{1,2,3\}$, let $v_h$ be a real number such that 
\[
x \in v_h + B \qqand (v_h +B) \cap (a_i+B) \neq \emptyset.   
\]
There exist elements $b_{p_h} \in B$ such that 
\[
v_h = x-b_{p_h} 
\]
and 
\[
(x-b_{p_h} +B ) \cap (a_i+B) \neq \emptyset. 
\]
It follows that there are elements $ b_{q_h}$ and $ b_{r_h}$ in $B$ such that 
\[
x-b_{p_h} + b_{q_h}  = a_i+b_{r_h}  
\] 
and so  
\[
x - a_i  = b_{p_h} - b_{q_h} + b_{r_h} 
\] 
for $h \in \{1,2,3\}$. 
We obtain  
\beq           \label{vonN:pq}
b_{p_1} + b_{q_2} + b_{r_1}  = b_{p_2} + b_{q_1} + b_{r_2} 
\eeq
and 
\beq           \label{vonN:pr}
b_{p_1} + b_{q_3} + b_{r_1}  = b_{p_3} + b_{q_1} + b_{r_3}. 
\eeq
Because the set $B$ is \Q-independent,  relation~\eqref{vonN:pq} implies 
\[
 b_{p_2} = b_{p_1} \quad \text{ or }  \quad    b_{p_2}= b_{q_2}  \quad  \text{ or }  \quad    b_{p_2} = b_{r_1}
\]
and  relation~\eqref{vonN:pr} implies  
\[
 b_{p_3} = b_{p_1} \quad \text{ or }  \quad    b_{p_3}= b_{q_3}  \quad  \text{ or }  \quad    b_{p_3} = b_{r_1}. 
\]

Suppose that $ b_{p_2} \neq b_{p_1}$.  If $b_{p_2} = b_{q_2} $, then 
\[
x-a_i = b_{r_2}
\]
and $x = a_i+b_{r_2} \in a_i + B$, which is absurd.  
Therefore, $ b_{p_2} \neq b_{p_1}$ implies $b_{p_2} = b_{r_1} $. 

Suppose that $ b_{p_3} \neq b_{p_1}$.  If $b_{p_3} = b_{q_3} $, then 
\[
x-a_i = b_{r_3}
\]
and $x = a_i+b_{r_3} \in a_i + B$, which is absurd.  
Therefore, $ b_{p_3} \neq b_{p_1}$ implies $b_{p_3} = b_{r_1} $. 
It follows that either $ b_{p_1} = b_{p_2}$ or $ b_{p_1} = b_{p_3}$ or 
$b_{p_2} = b_{r_1} = b_{p_3}$, 
and so at most two translates of $B$ contain $x$ and also intersect $a_i+B$ 
for some $i \in \{1,\ldots, k\}$.  Therefore, at most $2k$ translates of $B$ 
contain $x$ and intersect $\bigcup_{i=1}^k A_i$.  Because there are infinitely 
many translates of $B$ that contain $x$ and are contained in $J'\cap G$, 
the pigeon hole principle completes the proof of the Lemma. 
\end{proof}

\section{Finishing the proof of von Neumann's theorem} 
Let $J$ be a finite interval $J$ with endpoints $\inf(J) = c_0$ and $ \sup(J) = c_1$.   
Let 
\[
0 < 8 \varepsilon < \length(J) = c_1-c_0.
\]
Let $B = B(\varepsilon) = \{b_0, b_1, b_2, \ldots\}$ 
be the set and $G = G(\varepsilon)$ the group constructed in Lemma~\ref{vonN:lemma:Bset} 
with $\beta_0 = -2\varepsilon$ and $\beta_1 = 2\varepsilon$. 
The numbers $b_0 = \min(B)$ and $b_1 = \max(B)$ satisfy 
\[
2 \varepsilon < b_1 - b_0 < 4\varepsilon. 
\]
Let 
\[
C = J \cap G.
\]
Because $G$ is dense in $J$, we have $\inf(C) = c_0$ and $ \sup(C) = c_1$.
The set $C$ does not necessarily contain the numbers $c_0$ and $c_1$.  

If $c_0 \in C$, then  
\[
c_0 \in (c_0 - b_0) + B = \{c_0 - b_0 + b_n: b_n \in B\} \subseteq G. 
\] 
The inequality 
\[
c_0 \leq c_0 - b_0 + b_n \leq c_0 + (b_1 - b_0)  < c_0 + 4\varepsilon < c_1
\]
implies $c_0 - b_0 + b_n \in J$ for all $b_n \in B$ 
and so $ (c_0 - b_0) + B \subseteq J$. 
Therefore, 
\[
(c_0 - b_0) + B \subseteq J\cap G = C. 
\]
Similarly, if $c_1 \in C$, then 
\[
(c_1 - b_1) + B \subseteq J\cap G = C. 
\]

Suppose that $c_0 \in C$ and $c_1 \in C$. 
If 
\[
y \in \left( (c_0 - b_0) + B \right) \cap \left( (c_1 - b_1) + B  \right) 
\]
then there exist $b_p, b_q \in B$ such that 
\[
y = c_0 - b_0 + b_p = c_1 - b_1 + b_q 
\]
and so 
\[
8\varepsilon < c_1 - c_0 = (b_p-b_0)+ (b_1 - b_q) < 4\varepsilon + 4\varepsilon = 8\varepsilon
\]
which is absurd.  Therefore, if $c_0 \in C$ and $c_1 \in C$, then the translates 
$(c_0 - b_0) + B$ and $(c_1 - b_1) + B$ are disjoint. 

We shall now prove Theorem~\ref{vonN:theorem:JG=AB}.
Because $G$ is dense in $J$ and countably infinite, the subset $C$ of $G$ is countably infinite.  
Let 
\[
C = \{x_n:n = 0,1,2,\ldots \}
\]
be an enumeration of the elements of $C$ 
that satisfies the following initial conditions:  
\benum 
\item[(i)] If $c_0 \in C$ and $c_1 \in C$, 
then $x_0 = c_0$ and $x_1 = c_1$. 
Let $a_0 = x_0 - b_0$ and $a_1 = x_1-b_1$. 
We have $x_0 \in a_0+B$ and $x_1 \in a_1+B$.  
The translates $a_0+B$ and $a_1+B$ are disjoint and contained in $C$. 

\item[(ii)] If $c_0 \in C$ but $c_1 \notin C$, then $x_0 = c_0$.  
Let $a_0 = x_0 - b_0$.  We have $x_0 \in a_0+B \subseteq C$. 

\item[(iii)] If $c_0 \notin C$ 
but $c_1 \in C$, then $x_0 = c_1$.  
Let $a_0 = x_0-b_1$. We have $x_0 \in a_0+B \subseteq C$. 
\item[(iv)] 
If $c_0 \notin C$ and $c_1 \notin C$, then let $x_0$ be any element of $C$ 
and let $a_0 + B$ be any translate of $B$ such that $x_0 \in a_0+B \subseteq C$. 
Lemma~\ref{vonN:lemma:core} implies the existence of such translates. 
\eenum
We use induction to construct an infinite sequence $(a_j +B)_{j=0}^{\infty}$ 
of pairwise disjoint translates of $B$ such that $A\oplus B = C$, where $A = \{a_j: j=0,1,2,\ldots\}$. 

Suppose that $n \geq 0$ and that $\{ a_j + B: j=0,1,\ldots, n\}$ are pairwise disjoint translates 
of $B$ such that 
\[
\{x_0,x_1,\ldots, x_n\} \subseteq \bigcup_{j=0}^n (a_j +B) \subseteq C.  
\]
By Lemma~\ref{vonN:lemma:B-cover}, 
the set $ \bigcup_{j=0}^n (a_j +B)$ is a proper subset of $C$.
Let $\ell$ be the smallest positive integer such that 
$x_{\ell} \notin \bigcup_{j=0}^n (a_j +B) \subseteq C$. 
Note that $\ell \geq n+1$. 
By the core Lemma~\ref{vonN:lemma:core},  there is a translate $a_{n+1} + B \subseteq C$ 
such that 
\[
x_{\ell} \in a_{n+1} + B 
\]
\[
\left(  a_{n+1} + B \right) \cap \left(  \bigcup_{j=0}^n (a_j +B) \right) = \emptyset
\]
and 
\[
\left\{ x_0,x_1,\ldots, x_n, x_{n+1} \ldots, x_{\ell} \right\} \subseteq \bigcup_{j=0}^{n+1} (a_j +B) \subseteq C.
\]
Continuing inductively, we obtain a countably infinite set $A = \{a_0, a_1, a_2, \ldots \}$ 
such that 
\[
A+B = \bigcup_{j=0}^{\infty} (a_j+B) = C.
\]
This completes the proof of von Neumann's theorem.

\section{Extension to parallelepipeds}          \label{vonN:section:parallelepiped}

A parallelepiped in $\R^n$ is a set of the form
\[
P = \left\{\sum_{i=1}^n t_i \mbv_i: t_i \in I_i \text{ for all } i \in \{1,\ldots, n\} \right\}
\]
where $\{\mbv_1, \mbv_2,\ldots, \mbv_n\}$ is a linearly  independent set of vectors 
and $\{I_1, I_2, \ldots, I_n\}$ is a set of intervals that may be open, closed, or half open, half closed.  
We shall extend von Neumann's theorem to parallelepipeds.

\bt               \label{vonN:theorem:parallelepiped} 
Every parallelepiped in $\R^n$ is the union of countably infinitely many pairwise disjoint congruent sets.
\et

\begin{proof}
Let $\{ \mbe_1, \mbe_2,\ldots, \mbe_n\}$ be the standard basis in $\R^n$.
Let $\{I_1, I_2, \ldots, I_n\}$ be a set of intervals that may be open, closed, or half open, half closed.  
Consider the rectangular parallelepiped  
\[
P_0 = I_1 \times I_2 \times \cdots \times I_n 
= \left\{\sum_{i=1}^n t_i \mbe_i: t_i \in I_i \text{ for all } i \in \{1,\ldots, n\} \right\}  
\] 
and the parallelepiped  
\[
P = I_1 \times I_2 \times \cdots \times I_n 
= \left\{\sum_{i=1}^n t_i \mbv_i: t_i \in I_i \text{ for all } i \in \{1,\ldots, n\} \right\}  
\]
where $\{\mbv_1, \mbv_2,\ldots, \mbv_n\}$ is a linearly  independent set of vectors. 
By Theorem~\ref{vonN:theorem:Steinhaus}, 
there are infinite sets $A_1, A_2,\ldots, A_n$ 
and countably infinite sets $B_1, B_2,\ldots, B_n$ such that 
\[
I_i = A_i \oplus B_i
\]
for all $i \in \{1,\ldots, n \}$.  
Let 
\[
A = A_1 \times \cdots \times A_n \subseteq \R^n
\]
and 
\[
B = B_1 \times \cdots \times B_n \subseteq \R^n
\]
The set $B$ is countably infinite.  It follows from Lemma~\ref{vonN:lemma:k-sum}  that 
\[
P_0 = \prod_{i=1}^n I_i  = \prod_{i=1}^n ( A_i \oplus B_i ) 
= (A_1 \times \cdots \times A_n) \oplus (B_1 \times \cdots \times B_n)  = A \oplus B 
\]
and so $P_0$ is the union of countably infinitely many pairwise disjoint congruent sets.  

Define the linear isomorphism $T: \R^n \rightarrow \R^n$ by 
$T(\mbe_i) = \mbv_i$ for all $i \in \{1,2,\ldots, n\}$. 
Then $T(P_0) = P$.  By Lemma~\ref{vonN:lemma:Tisom}, 
the parallelepiped  $P$ is the union of countably infinitely many pairwise disjoint congruent sets. 
This completes the proof. 
\end{proof}

It is not clear if decomposition into countably infinitely many pairwise disjoint congruent sets 
is possible for other geometrical objects, such as, for example, triangles in the plane.


\end{document}